\documentclass[11pt]{article}

\usepackage{graphicx}
\usepackage{amssymb}
\newtheorem{definition}{Definition}

\newtheorem{theorem}{Theorem}
\newtheorem{lemma}{Lemma}

\newcommand{\EQ}{\begin{equation}\begin{array}{lllllllll}}
\newcommand{\EE}{\end{array}\end{equation}}
\newcommand{\EQQ}{$$\begin{array}{lllllllll}}
\newcommand{\EEE}{\end{array}$$}
\newcommand{\MT}{\left[ \begin{array}{ccccccccc}}
\newcommand{\EM}{\end{array}\right]}
\newcommand{\bean}{\begin{equation}\begin{array}{rcllllllll}}
\newcommand{\eean}{\end{array}\end{equation}}
\newcommand{\bea}{$$\begin{array}{rcllllllll}}
\newcommand{\eea}{\end{array}$$}
\def\ds{\displaystyle}
\def\Fr{\ds \frac}
\def\=={&=&}
\def\dsum{\ds\sum}
\def\Real{{I\!\!R}}
\def\calL{{\cal L}}
\def\calH{{\cal H}}
\def\calB{{\cal B}}
\def\L2O{X}
\def\DBL{D_B(\calL)}
\def\DB{D_B}

\title{\LARGE \bf
The Consistency of Partial Observability for PDEs
\thanks{This work was supported in part by NRL and AFOSR}
}

\author{Wei Kang \thanks{Wei Kang is with Faculty of Applied Mathematics, Naval Postgraduate School, Monterey, CA, USA
        {\tt\small wkang@nps.edu}} 
}

\begin{document}

\maketitle
\thispagestyle{empty}

\begin{abstract}
In this paper, a new definition of observability is introduced for PDEs. It is a quantitative measure of partial observability. The quantity is proved to be consistent if approximated using well posed approximation schemes. A first order approximation of an unobservability index using empirical gramian is introduced. For linear systems with full state observability, the empirical gramian is equivalent to the observability gramian in control theory. The consistency of the defined observability is exemplified using a Burgers' equation. 

\end{abstract}

\section{Introduction}
Observability is a fundamental property of dynamical systems \cite{isidori,kailath} that has an extensive literature. It is a measure of well-posedness for the estimation of system states using both sensor information and non-sensor knowledge. It is impossible in this paper to review all the main results on this subject while some interesting work can be found in \cite{gauthier,hermann,xia} for nonlinear systems, \cite{zuazua} for PDEs, \cite{mohler} for stochastic systems, and \cite{zheng} for normal forms. For complicated problems, the challenge is to define the concept so that it captures the fundamental nature of observability, and meanwhile, it is practically verifiable. In \cite{kang-xu1,kang-xu2}, a definition of observability is introduced using dynamic optimization as a tool. This concept  is developed in a project of finding the best sensor locations for data assimilation, a computational algorithm widely used in numerical weather prediction. Different from traditional definitions of observability, the one in \cite{kang-xu1,kang-xu2} is able to collectively address several issues in a unified framework, including quantitatively measure observability, determine partial observability, accommodate both sensor and non-sensor information, etc. Moreover, computational methods of dynamic optimization provide practical tools to numerically approximate the observability of complicated systems that cannot be treated using an analytic approach. 

To extend the defintion of observability in \cite{kang-xu1,kang-xu2} to systems defined by PDEs, several fundamental issues must be addressed. A partial observability makes perfect sense for infinite dimensional systems such as PDEs. However, its computation must be carried out using finite dimensional approximations. It is known in the literature that an ODE approximation of a PDE may not preserve the property of observability, even if the approximation scheme is convergent and stable \cite{zuazua,cohn}. Therefore, to apply the concept of observability to PDEs, it is important to understand its consistency in ODE approximations. 
 
In Section \ref{sec2}, some examples from existing literature are introduced to illustrate the issues being addressed in this paper. Observability is defined for PDEs in Section \ref{sec3}. In Section \ref{sec4}, a theorem on the consistency of observability is proved. The relationship between the unobservability index and an empirical gramian is addressed in Section \ref{sec5}, which serves as a first order approximation of the observability. In Section \ref{sec6}, the concept is extended to nonlinear systems and a theorem of consistency is proved. The theory is verified using a Burgers' equation in Section \ref{sec7}.

\section{Some issues on observability}
\label{sec2}
Consider the initial value problem of a heat equation
$$\begin{array}{lll}
u_t(x,t)=u_{xx}(x,t), \; x\in [0, L], t\in [0, T]\\
u(0,t)=u(L, t)=0\\
u(x,0)=f(x)
\end{array}$$
Suppose the measured output is 
$$y(t)=u(x_0,t)$$
for some $x_0\in [0,L]$. In this example, we assume $L=2\pi$, $T=10$, and $x_0=0.5$. Suppose the solution and its output has the following form
$$\begin{array}{lll}
u(x,t)=\ds\sum_{k=1}^{\infty} \bar u_k(t)\sin \left(\Fr{k\pi x}{L}\right)\\
y(t)=\ds\sum_{k=1}^{\infty} \bar u_k(t)\sin \left(\Fr{k\pi x_0}{L}\right)
\end{array}$$
Then the Fourier coefficients satisfy the following ODE
$$\Fr{ d\bar u_k}{dt}=-\left(\Fr{k\pi}{L}\right)^2 \bar u_k$$
Define 
$$\begin{array}{lll}
u^N=\MT \bar u_1 & \bar u_2 & \cdots \bar u_N\EM^T\\
A^N=\mbox{diag}\left(\MT \left(\Fr{\pi}{L}\right)^2 & \left(\Fr{2\pi}{L}\right)^2&\cdots &\left(\Fr{N\pi}{L}\right)^2\EM\right)\\
C^N=\MT \sin\left(\Fr{\pi x_0}{L}\right) & \sin\left(\Fr{2\pi x_0}{L}\right)&\cdots &\sin\left(\Fr{N\pi x_0}{L}\right)\EM
\end{array} $$
A $N$th order approximation of the original initial value problem with output is defined by a system of ODEs
\EQ
\label{eqheatode}
\Fr{du^N}{dt}=-A^Nu^N\\
u^N(0)=u^N_0\\
y=C^Nu^N
\EE
A gramian matrix \cite{kailath} can be used to measure the observability of $u^N_0$. More specifically, given $N>0$ the observability gramian is
$$W=\ds\int_0^Te^{(-A^N)'t}(C^N)'C^Ne^{-A^Nt}dt$$
Its smallest eigenvalue, $\sigma^N_{min}$, measures the observability of $u^N_0$. A small value of $\sigma^N_{min}$ implies weak observability. If the maximum sensor error is $\epsilon$, then the worst estimation error of $u^N_0$ is 
\EQ
\label{esterror}
\Fr{\epsilon}{\sqrt{\sigma^N_{min}}}
\EE

\begin{figure}[!ht]
	\begin{center}
		\includegraphics[width=3.0in,height=2.0in]{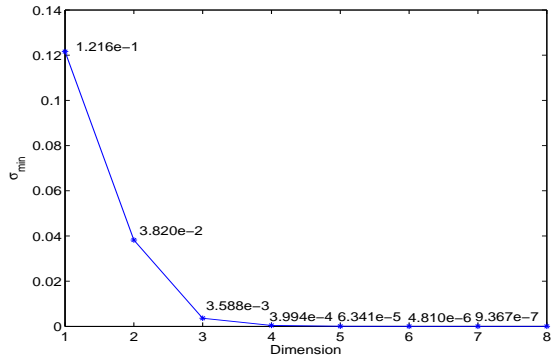} 
		\caption{Observability of heat equation}
		\label{Figheat}
		\end{center}
\end{figure}

The system has infinitely many modes in its Fourier expansion. However, it has a single output. The computation shows that the output can make the first mode observable. However, when the number of modes is increased, their observability decreases rapidly. From Figure \ref{Figheat}, for $N=1$ we have $\sigma^N_{min}=1.216$, which implies a reasonably observable $\bar u_1(0)$. However, when $N$ is increased, the observability decreases rapidly. For $N=8$, $\sigma^N_{min}$ is almost zero, i.e 
$$\MT \bar u_1(0)&\bar u_2(0)&\cdots &\bar u_8(0)\EM^T$$ 
is extremely weakly observable, or practically unobservable. According to (\ref{esterror}), a small sensor noise results in a big estimation error.  

The family of solutions of a PDE is an infinite dimensional space. Finite many sensors may not provide adequate information to accurately estimate all modes in an initial state. In the example of heat equation, a single output makes the first two or three modes observable. All other modes are practically unobservable. On the other hand, in practical applications a finite number of modes is enough to provide accurate approximations. All we need is the observability of the finite modes. This is the reason we would like to {\it measure} a system's {\it partial observability}. The concept is useful not only for PDEs. In large scale systems with high (finite) dimensions, it may not be possible or necessary to make the entire state space observable. A partial observability is all we need.  

Another issue to be addressed in this paper is {\it consistency}. In general the observability for PDEs is numerically computed using a system of ODEs as an approximation. However, it is not automatically guaranteed that the observability of the ODEs is consistent with the observability of the original PDE. In fact, a convergent discretization of a PDE may not preserve its observability. Take the following wave equation as an example
\EQ
\label{wave}
u_{tt}-u_{xx}=0, &0<x<L, \; 0<t<T\\
u(0,t)=u(L,t)=0, & 0<t<T\\
u(x,0)=u_0(x), u_t(x,0)=u_1(x), &0<x<L
\EE
It is known that the total energy of the system can be estimated by using the energy concentrated on the boundary. However, in \cite{zuazua} it is proved that the discretized ODEs do not have the same observability. The energy of solutions is given by
$$E(t)=\Fr{1}{2}\ds\int_0^L \left( |u_t(x,t)|^2+|u_x(x,t)|^2  \right) dx$$
This quantity is conserved along time. It is know that, when $T>2L$, the total energy of solutions
can be estimated uniformly by means of the energy concentrated on the boundary $x=L$. More precisely, there exists $C(T)>0$ such that 
\EQ
\label{waveobs}
E(0) \leq C(T)\ds\int_0^T |u_x(L,t)|^2dt
\EE
Now consider the discretized system using a finite difference method,
\EQ
\label{wave2}
u''_j(t)=\Fr{u_{j+1}(t)+u_{j-1}(t)-2u_j(t)}{h^2}, & 0<t<T, \;j=1,2,\cdots,N\\
u_0(t)=u_{N+1}(t)=0, & 0<t<T\\
u_j(0)=u^0_j, u'_j(0)=u^1_j, &j=0,1,\cdots,N+1
\EE
The total energy of the ODEs is given by
$$E_h(t)=\Fr{h}{2}\ds\sum_{j=0}^N \left(|u'_j(t)|^2+\left|\Fr{u_{j+1}(t)-u_j(t)}{h}   \right|^2   \right)$$
This quantity is conserved along the trajectories of the ODEs. The energy on the boundary is defined by
$$\ds\int_0^T \left| \Fr{u_N(t)}{h} \right|^2dt$$
Because the solutions of (\ref{wave2}) converges to the solutions of (\ref{wave}), we expect that the total energy of (\ref{wave2}) can be uniformly estimated using the energy concentrated along a boundary, i.e. the following inequality similar to (\ref{waveobs}) holds for some $C(T)$ uniformly in $h$,
$$E_h(0)\leq C(T)\ds\int_0^T \left| \Fr{u_N(t)}{h} \right|^2dt$$
However, it is proved in \cite{zuazua} that this inequality is not true. In fact, the ratio between the total energy and the energy along the boundary is unbounded as $h\rightarrow 0$. To summarize, the observability of a PDE is not necessarily preserved in its discretizations. 

In this paper, we introduce a {\it quantitative measure} of {\it partial observability} for PDEs. Sufficient conditions are proved for the {\it consistency} of the observability for well posed discretization schemes. 

\section{Problem formulation}
\label{sec3}
\setcounter{equation}{0}
Following \cite{canuto}, we formulate a linear evolution problem
\EQ
\label{eqpdemodel}
u_t+\calL u=g, &\mbox{in } \Omega \times (t_0, t_1]\\
\calB u=0 &\mbox{on } \partial \Omega_b \times (t_0, t_1]\\
u=u_0 & \mbox{in } \Omega \times \{ t=0\}\\
y(t)=\calH (u(\cdot,t))
\EE
where $\Omega$ is an open set in $\Real^n$, $\calL$ is a linear operator, bounded or unbounded, in its domain $ D(\calL)$, which is a subspace of a Banach space $(X, ||\cdot||_X)$. In the following, $u(t)$ represents $u(\cdot,t)\in D(\calL)$. The boundary condition is defined by a set of linear boundary operators on a part of $\partial \Omega$, denoted by $\partial \Omega_b$. To take into consideration of the boundary condition, we define
$$D_B(\calL)=\{v\in D(\calL)|\; \calB v =0 \mbox{ on } \partial \Omega_b\}$$
We consider the family of solutions, in either strict or weak sense, generated by $u(0)\in D_{B}(\calL)$. The right-hand side $g$ is a continuous function of the variable $t$ with values in $\L2O$, i.e $g\in C([t_0,t_1],\L2O)$. A solution $u(t)$ for this problem is a $\L2O$ valued function that is continuous in $[t_0,t_1]$, $du/dt$ exists and is continuous for all $t\in (t_0,t_1]$, satisfying $u(0)=u_0$, and $u(t)\in \DBL$ for all $t\in (0,t_1]$. We assume that (\ref{eqpdemodel}) is a {\it well posed problem} in the Hadamard sense (\cite{hille,richtmyer}). More specifically, 
\begin{itemize}
\item For any $u_0\in D_{B}(\calL)$, (\ref{eqpdemodel}) has a solution.
\item The solution for $u_0\in D_{B}(\calL)$ is unique.
\item The solution of (\ref{eqpdemodel}) depends continuously on its initial value.
\end{itemize}
In (\ref{eqpdemodel}), $y(t)=\calH (u(\cdot,t))$ represents the output of the system in which $\calH$ is a linear operator from $\L2O$ to $\Real^p$. The output $y(t)\in C([t_0,t_1])$ has a norm denoted by $|| y||_Y$. Rather than the entire state space, the observability is defined in a finite dimensional subspace. Let 
$$W=\mbox{span}\{ e_1,e_2,\cdots,e_s\}$$
be a subspace of $D_{B}(\calL)$ generated by a basis $\{ e_1,e_2,\cdots,e_s\}$. In the following, we analyze the observability of $u(0)$ using estimates from $W$. Therefore $W$ is called the {\it space for estimation}.

Let $u(t)$ be a solution of (\ref{eqpdemodel}). Suppose $u_w(0)$ generate the best estimation of $u(0)$ in $W$ in the sense that $u_w(t)$ minimizes the following output error,
\EQ
\label{uproj}
\min||\calH (u_w)-\calH (u)||_Y\\
\mbox{subject to} \\
du_{w}/dt+\calL u_w =g\\
u_w (t)\in \DBL\mbox{ for all } t\in (t_0,t_1]\\
u_w(0)\in W
\EE
Let $u_r(t)=u(t)-u_w(t)$ be the remainder, then 
\EQ
\label{decomposition}
u(t)=u_w(t)+u_r(t)
\EE
According to (\ref{uproj}), the best estimate of $u(0)$ in $W$ is $u_w(0)$. The estimation error, $||u_r(0)||_X$, is not directly caused by the output noise, i.e. this error of remainder cannot be reduced no matter how accurate the output is measured. Therefore, the following {\it partial observability} is defined for $u_w(0)$ only. Or equivalently, we assume $u(0)\in W$ in the definition. 

\begin{definition}
\label{def1}
Given a nominal trajectory $u$ of (\ref{eqpdemodel}) with $u(0)\in W$. For a given $\rho >0$, define
\EQ
\label{eqcost}
\epsilon = \inf ||\calH (\hat u)-\calH (u)||_Y
\EE
where $\hat u$ satisfies
\EQ
\label{hatu}
\hat u_t+\calL \hat u =g\\
\hat u (t)\in \DBL\mbox{ for all } t\in (t_0,t_1]\\
\hat u(0)\in W\\
||\hat u(0) -u(0)||_X=\rho 
\EE
Then $\rho/\epsilon$ is called the unobservability index of $u(0)$ along the trajectory $u(t)$.
\end{definition}

\noindent{\bf Remark}. The ratio $\rho/\epsilon$ can be interpreted as follows: if the maximum error of the measured output, or sensor error, is $\epsilon$, then the worst estimation error of $u(0)$ is $\rho$. Therefore, a small value of $\rho/\epsilon$ implies strong observability of $u(0)$. \\
\\
{\bf Remark}. If $u(0)$ is not in $W$, then the overall error in the estimation of $u(0)$ is bounded by the error of the estimate of $u_w$ plus the remainder, or the truncation error, $||u_r(0)||_X$. For $u(0)$ to be strongly observable, it requires a strong observability of $u_w(0)$ and a small remainder $u_r(0)$.\\
\\
{\bf Remark}. If $u(0)=r_w(0)+u_r(0)$ is not in $W$, in order to compute the observability of $u_w(0)$ one must solve (\ref{uproj}) first to find $u_w(0)$. In fact, this is not necessary. Given any solution of (\ref{hatu}). It can be expressed as
$$\hat u(t)=u_w(t)+\delta u(t)$$
where $\delta u(t)$ is a solution of the associated homogeneous PDE and
$$\delta u(0) \in W, \;\; ||\delta u(0)||_X=\rho$$
Therefore, 
$$\begin{array}{lll}
\hat u(t)+u_r(t)=u_w(t)+u_r(t)+\delta u(t)\\
=u(t)+\delta u(t)
\end{array}$$
Because 
$$\begin{array}{lll}||\calH(\hat u)-\calH(u_w)||_Y\\
=||\calH(\hat u+u_r)-\calH(u_w+u_r)||_Y\\
=\calH( u+\delta u)-\calH(u)||_Y
\end{array}
$$
the dynamic optimization (\ref{eqcost})-(\ref{hatu}) is equivalent to 
\EQ
\label{eqcost2}
\epsilon = \inf ||\calH (\hat u)-\calH (u)||_Y\\
\mbox{subject to}\\
\hat u_t+\calL \hat u =g\\
\hat u (t)\in \DBL\mbox{ for all } t\in (t_0,t_1]\\
\hat u(0)\in u+W\\
||\hat u(0) -u(0)||_X=\rho 
\EE
In (\ref{eqcost2}), it is not necessary to compute $u_w(t)$. \\
\\
{\bf Remark}. It can be shown that, for linear problems, $\rho/\epsilon$ is a constant. The expressions in (\ref{eqcost})-(\ref{hatu}) can be simplified (See Section \ref{sec5}). However, we prefer the form adopted in Definition \ref{def1} because it can be easily generalized to nonlinear problems\\

To numerically compute a system's observability, (\ref{eqpdemodel}) is approximated by ODEs. In this paper, we consider a general approximation scheme using a sequence of ODEs,
\EQ
\label{eqodemodel}
du^N/dt+A^Nu^N=g^N, & u^N\in \Real^N\\
u^N(0)=u^N_0
\EE
where $N\geq N_0$ for some integer $N_0$. The approximation is constructed using two linear mappings
\EQ
\label{mapping}
P^N: D(\calL) \rightarrow \Real^N\\
\Phi^N: \Real^N \rightarrow \L2O
\EE
In addition, a norm, $||\cdot ||_N$, is defined on $\Real^N$. 
The approximation scheme is said to be {\it well posed} if it satisfies the following conditions.
\begin{itemize}
\item Given any solution $u(t): [t_0,t_1]\rightarrow \DBL$, of (\ref{eqpdemodel}). Suppose $u^N(t)$ is a solution of (\ref{eqodemodel}) satisfying $u^N(0)=P^N(u(0))$, then
\EQ
\label{mapping2}
\ds\lim_{N\rightarrow \infty} ||\Phi^N(u^N(t))-u(t)||_X =0
\EE
converges uniformly on $[t_0,t_1]$.
\item For any $u\in \L2O$, the sequence defined by $u^N=P^Nu$ satisfies
\EQ
\label{mapping3}
\ds\lim_{N\rightarrow \infty}||u^N||_N=||u||_X
\EE
\end{itemize}

Given the space for estimation $W$, we define a sequence of subspaces, $W^N \subseteq\Real^N$, by
$$W^N=P^N(W)$$
They are used as the space for estimation in $\Real^N$. If $\{e_1,e_2,\cdots,e_s\}$ is a basis of $W$, then their projections to $W^N$ are denoted by
$$e_i^N=P^N(e_i), \;\; i=1,2,\cdots,s$$  
So $W^N=\mbox{span}\{ e_1^N, e_2^N,\cdots,e_s^N\}$.\\

\noindent{\bf Example}.
For a spectral method, approximate solutions can be expressed in terms of an orthonormal basis 
$$\{ q_N(x): N=0,1,2,\cdots \}$$
For any function, 
$$v(x)=\sum_{k=0}^\infty v_kq_N(x)\in D(\calL)$$ one can define 
$$P^N (v)=\MT v_0&v_1&\cdots&v_N\EM^T$$
Obviously, $\Phi^N$ is defined by
$$\Phi^N(\MT v_0&v_1&\cdots&v_N\EM^T)=\sum_{k=0}^N v_kq_N$$
Because the basis is orthonormal, the $l^2$ norm and the inner product in $\Real^N$ is consistent with those in $L^2(\Omega)$. 
$\square$\\

\noindent{\bf Example}. Some approximation methods, such as finite difference and finite element, are based on a grid defined by a set of points in space, 
$$\{ x_k\}_{k=1}^N$$
and a basis $\{ q_k \}$ satisfying
$$q_k(x_j)=\left\{ \begin{array}{lll} 1 & k=j\\
0, & otherwise
\end{array}\right.$$
In this case, the mappings in the approximation scheme is defined as follows
$$\begin{array}{lll}
P^N(v)=\MT v(x_1)&v(x_2)&\cdots&v(x_N)\EM^T\\
\Phi^N(\MT v_1&v_2&\cdots&v_N\EM^T)=\ds\sum_{k=1}^N v_kq_k
\end{array}$$
The inner product in $\Real^N$ can be induced from the $L^2$ space, i.e. for $u, v\in \Real^N$, 
$$<u,v>_N=<\ds\sum_{k=1}^N u_kq_k, \ds\sum_{k=1}^N v_kq_k>$$
If the functions in $D(\calL)$ are uniformly continuous, then $\Phi^N(P^N(v))$ converges to $v$ uniformly. It can be shown that $<\cdot, \cdot>_N$ is consistent with the inner product in $L^2(\Omega)$.$\square$\\

Following \cite{kang-xu1}, we define the observability for ODE systems. 

\begin{definition}
\label{def2}
Given $\rho >0$ and a trajectory $u^N(t)$ of (\ref{eqodemodel}) with $u^N(0)\in W^N$. Let 
$$\epsilon^N = \inf ||\calH\circ\Phi^N(\hat u^N)-\calH\circ \Phi^N(u^N)||_Y$$
where $\hat u^N$ satisfies
\EQ
\label{hatuN}
d\hat u^N/dt+A^N\hat u^N=g^N\\
\hat u^N(0)\in W^N\\
||\hat u^N(0) -u^N(0)||_N=\rho 
\EE
Then $\rho/\epsilon^N$ is called the unobservability index of $u^N(0)$ in the space $W^N$.
\end{definition}

\section{The consistency of observability}
\label{sec4}
\setcounter{equation}{0}

In this section, we address the consistency of observability. The theorem is based on a continuity assumption. In the problem formulation, the output mapping $\calH$ is not required to be bounded. However, in the proof of the consistency theorem we require $\calH$ be continuous in the following subspace of $\L2O$ extended from $W$
$$W_E=\mbox{span}\left\{ \{ e_1,e_2,\cdots,e_s\}\cup \{ \Phi^N(e^N_1),\cdots,\Phi^N(e^N_s)\}_{N=N_0}^\infty\right\} $$

\noindent{\bf Output Continuity Assumption}: Given a sequence 
$$\{v_k(t)\}_{k=k_0}^\infty \subset C^1([t_0,t_1],W_E)$$ 
If $v_k(t)$ converges to $v(t)$ in $W_E$ uniformly on $[t_0,t_1]$, then 
$$\ds\lim_{k\rightarrow \infty}\calH(v_k)=\calH(v)$$ 

This assumption is automatically satisfied if $\calH$ is a bounded linear operator. Or if $W_E$ has a finite dimension, then this assumption is satisfied even if the operator is unbounded in $\L2O$.

\begin{theorem} 
\label{theorem1}
Suppose the initial value problem (\ref{eqpdemodel}) and its approximation scheme (\ref{eqodemodel})-(\ref{mapping}) are well posed. Suppose $\calH$ satisfies Output Continuity Assumption. Then, 
\EQ
\label{consist}
\ds\lim_{N\rightarrow \infty} \epsilon^N =\epsilon
\EE
\end{theorem}

To prove this theorem, we need the following lemma.
\begin{lemma}
\label{lemma1}
Given a sequence $\hat u^{N}(t)$, $N\geq N_0$, satisfying (\ref{hatuN}). Then there exists a subsequence, $\hat u^{N_k}(t)$, so that $\{ \Phi^{N_k}(\hat u^{N_k}(t))\}_{k=1}^\infty$ converges uniformly to a solution of (\ref{hatu}). 
\end{lemma}

{\it Proof }. Let $\bar u(t)$ be a solution of the original PDE with initial value
$$\bar u(0)=0$$
and let $h_i(t)$, $i=1,2,\cdots,s$, be the solution of the associated homogeneous PDE with initial value $e_i$, i.e.
\EQ
\label{eqpdehomo}
\partial h_i/\partial t+\calL h=0, &\mbox{in } \Omega \times (t_0, t_1]\\
\calB u=0 &\mbox{on } \partial \Omega_b \times (t_0, t_1]\\
u=e_i & \mbox{in } \Omega \times \{ t=0\}
\EE
Then any solution of the nonhomogeneous PDE with an initial value in $W$ has the form 
\EQ
\label{eqsuperp}
\bar u + \ds\sum_{i=1}^s a_ih_i
\EE
For each $N$, define $\bar u^N(t)$ be the solution of an approximating ODE satisfying
$$\bar u^N(0)=0$$
We know that $\Phi^N(\bar u^N)$ approaches $\bar u$ uniformly as $N\rightarrow \infty$.
Let $h_i^N(t)$, $i=1,2,\cdots,s$, be the solution of the associated homogeneous ODE with initial value $e^N_i$, i.e.
\EQ
\label{eqodehomo}
dh_i^N/dt+A^N h_i^N=0\\
h_i^N(0)=e^N_i
\EE
Then $\Phi^N(h^N_i(t))$ approaches $h_i(t)$ uniformly as $N\rightarrow \infty$. For each $\hat u^N(t)$ in Lemma \ref{lemma1}, it can be expressed as
$$\hat u^N(t)=\bar u^N(t)+\ds\sum_{i=1}^s a_i^N h_i^N(t)$$
and 
$$\hat u^N(0)=\ds\sum_{i=1}^s a_i^N e_i^N$$
From the initial condition in (\ref{hatuN}), we know that the set 
$$\{||\ds\sum_{i=1}^s a_i^N e_i^N||_N, \;\; N\geq N_0\}$$
is bounded. Because the norms are consistent, we can prove that
the sequence $\{(a_1^N,a_2^N,\cdots,a_s^N)^T\}_{N=N_0}^\infty$ has a bounded subsequence which converges under the standard norm $||\cdot||_2$. Let $\{(a_1^{N_k},a_2^{N_k},\cdots,a_s^{N_k})^T\}_{k=1}^\infty$ be the convergent subsequence with a limit $(a_1,a_2,\cdots,a_s)^T$. Then
$$\begin{array}{rcl}
&&\ds\lim_{k\rightarrow \infty} \Phi^{N_k}(\hat u^{N_k}(t))\\
&=&\ds\lim_{k\rightarrow \infty} (\Phi^{N_k}(\bar u^{N_k}(t))+\ds\sum_{i=1}^s a_i^{N_k} \Phi^{N_k}(h_i^{N_k}(t)))\\
&=&\bar u(t)+\ds\sum_{i=1}^s a_ih_i(t)\\
&\triangleq& \hat u(t)
\end{array}$$
The limit converges uniformly. From (\ref{eqsuperp}), $\hat u(t)$  must be a solution of the PDE in (\ref{hatu}). Because 
$$\begin{array}{lll}
||\hat u^{N_k}(0)-u^{N_k}(0)||_{N_k}
=\rho
\end{array}$$ 
and because of the consistency of the norms, we have 
$$\begin{array}{lll}
|| \hat u(0)-u(0)||_X\\
=\ds\lim_{k\rightarrow \infty} || \ds\sum_{i=1}^s a_i e_i^{N_k} -u^{N_k}(0)||_{N_k}\\
=\ds\lim_{k\rightarrow \infty} || \ds\sum_{i=1}^s a_i^{N_k} e_i^{N_k} -u^{N_k}(0)||_{N_k}\\
=\rho
\end{array}$$ 
Therefore, $\{ \Phi^{N_k}(\hat u^{N_k})\}_{k=1}^\infty$ 
converges uniformly to a solution of (\ref{hatu}).
$\square$\\

{\it Proof of Theorem \ref{theorem1}}. First, we prove
\EQ
\label{lowlim}
\liminf \epsilon^N \geq \epsilon
\EE
Suppose this is not true, then $\liminf \epsilon^N < \epsilon$. There exists $\alpha >0$ and a subsequence $N_k\rightarrow \infty$ so that
$$\epsilon^{N_k} < \epsilon -\alpha$$
for all $N_k$. From the definition of $\epsilon^N$, there exist $\hat u^{N_k}(t)$ satisfying (\ref{hatuN}) such that 
$$||\calH\circ\Phi^{N_k}(\hat u^{N_k})-\calH\circ\Phi^{N_k}(u^{N_k})||_Y < \epsilon - \alpha$$
From Lemma \ref{lemma1}, we can assume that $\Phi^{N_k}(u^{N_k})$ converges to $\hat u$ uniformly and $\hat u$ satisfies (\ref{hatu}). From Output Continuity Assumption, 
$$\lim_{k\rightarrow \infty} ||\calH\circ\Phi^{N_k}(\hat u^{N_k})-\calH\circ\Phi^{N_k}(u^{N_k})||_Y =||\calH(\hat u)-\calH(u)||_Y \leq \epsilon - \alpha$$
However, from the definition of $\epsilon$, we know 
$$ \epsilon \leq ||\calH(\hat u)-\calH(u)||_Y$$
A contradiction is found. Therefore, (\ref{lowlim}) must hold. 

In the next step, we prove
\EQ
\label{uplim}
\limsup \epsilon^N \leq \epsilon
\EE
It is adequate to prove the following statement: for any $\alpha>0$, there exists $N_1>0$ so that 
\EQ
\label{ineq}
\epsilon^N < \epsilon +\alpha
\EE
for all $N\geq N_1$. From the definition of $\epsilon$, there exists $\hat u$ satisfying (\ref{hatu}) so that 
\EQ
\label{eqad1}
||\calH(\hat u)-\calH(u)||_Y < \epsilon +\alpha
\EE
Let $\hat u^N$ be a solution of the ODE
$$d \hat u^N/dt+A^N\hat u^N=g^N$$
with the initial value
$$\hat u^N(0)=P^N(\hat u(0))$$ 
Then the following limit converges uniformly
\EQ
\label{eqlim}
\ds\lim_{N\rightarrow \infty} ||\Phi^N(\hat u^N(t))-\hat u(t)||_X = 0
\EE
A problem with $\hat u^N(0)$ is that its distance to $u^N(0)$ may not be $\rho$, which is required by (\ref{hatuN}). Let us define 
$$\begin{array}{lll}
\bar u^N(t)=\gamma_N(\hat u^N(t)-u^N(t))+u^N(t)\\
\gamma_N=\Fr{\rho}{||\hat u^N(0)-u^N(0)||_N}
\end{array}$$
Then $\bar u(t)$ satisfies $(\ref{hatuN})$. Due to the consistency of the norms and the fact $||\hat u(0)-u(0)||_X=\rho$, we know
$$\lim_{N\rightarrow \infty} \gamma_N=1$$
From (\ref{eqlim}), we have
$$\Phi^N(\bar u^N(t))-\Phi^N(u^N(t))=\gamma_N(\Phi^N(\hat u^N(t))-\Phi^N(u^N(t)))$$
converges uniformly to $\hat u-u$. Output Continuity Assumption and (\ref{eqad1}) imply
$$\begin{array}{lll}
\ds\lim_{N\rightarrow \infty} ||\calH\circ\Phi^N(\bar u^N)-\calH\circ\Phi^N(u^N)||_Y\\
=||\calH(\hat u)-\calH (u)||_Y\\
<\epsilon + \alpha
\end{array}$$
This implies that there exits $N_1>0$ so that
$$||\calH\circ\Phi^N(\bar u^N)-\calH\circ\Phi^N(u^N)||_Y < \epsilon+\alpha$$
for all $N\geq N_1$. From the definition of $\epsilon^N$, we know
$$\epsilon^N \leq ||\calH\circ\Phi^N(\bar u^N)-\calH\circ\Phi^N(u^N)||_Y< \epsilon + \alpha$$
for all $N>N_1$. Therefore, (\ref{uplim}) holds. 

To summarize, the inequalities (\ref{lowlim}) and (\ref{uplim}) imply 
$$\lim_{N\rightarrow \infty} \epsilon^N =\epsilon$$
$\square$

\section{Empirical gramian matrix}
\label{sec5}
\setcounter{equation}{0}

In this section, we assume that $W^N$ and the space of $y(t)$ are both Hilbert spaces with inner products $<,>_N$ and $<\cdot,\cdot>_Y$, respectively. Then $\epsilon^N/\rho$ equals the smallest eigenvalue of a gramian matrix. More specifically, let 
$$\{e_1^N,e_2^N,\cdots,e_s^N\}$$
be a set of orthonormal basis of $W^N$. For any $\hat u^N$ satisfying (\ref{hatuN}), we have
$$\begin{array}{lll}
\hat u^N(t)-u^N(t)\\
=e^{-A^Nt}(\hat u^N(0)-u^N(0))\\
=\ds\sum_{k=1}^s a_k e^{-A^Nt}e_k^N
\end{array}$$
for some coefficients satisfying
$$\ds\sum_{k=1}^s a_k^2=\rho$$
Therefore, 
$$\begin{array}{lll}
<\calH\circ\Phi^N(\hat u^N (t)-u^N(t)),\calH\circ\Phi^N(\hat u^N (t)-u^N(t))>_N\\
=<\ds\sum_{k=1}^s a_k\calH\circ\Phi^N(e^{-A^Nt}e_k^N),\ds\sum_{l=1}^s a_l\calH\circ\Phi^N(e^{-A^Nt}e_l^N)>_N\\
=\MT a_1&a_2&\cdots&a_s\EM G \MT a_1&a_2&\cdots&a_s\EM^T
\end{array}$$
where $G$ is the gramian 
\EQ
\label{eqgramian}
G=\MT G_{ij} \EM_{s\times s}, & G_{ij}=<\calH\circ\Phi^N(e^{-A^Nt}e^N_i),\calH\circ\Phi^N(e^{-A^Nt}e^N_j)>_Y
\EE
This matrix is the same as the observability gramian if $W^N$ is the entire space and if $y(t)$ lies in a $L^2$ space. It is straightforward to prove
$$\begin{array}{rclll}
(\epsilon^N)^2
&=&\min <\calH\circ\Phi(\hat u (t)-u(t)),\calH\circ\Phi(\hat u (t)-u(t))>_N\\
&=&\ds\min_{\sum a_k^2=\rho^2} \MT a_1&a_2&\cdots&a_s\EM G \MT a_1&a_2&\cdots&a_s\EM^T\\
&=&\sigma_{min}\rho^2
\end{array}$$
where $\sigma_{min}$ is the smallest eigenvalue of the gramian $G$. To summarize, if $u^N(0)$ and $y(t)$ lie in Hilbert spaces, then the unobservability index of the discretized system can be computed using the smallest eigenvalue of the gramian (\ref{eqgramian})
\EQ
\label{unobsgramian}
\rho/\epsilon^N=\Fr{1}{\sqrt{\sigma_{min}}}
\EE
If $\{e_1^N,e_2^N,\cdots,e_s^N\}$ is not an orthonormal basis, then 
$$\begin{array}{rclll}
(\epsilon^N)^2
&=&\min \MT a_1&a_2&\cdots&a_s\EM G \MT a_1&a_2&\cdots&a_s\EM^T\\
\mbox{subject to}\\
 \rho^2&=&\MT a_1&a_2&\cdots&a_s\EM S \MT a_1&a_2&\cdots&a_s\EM^T
\end{array}$$
where 
\EQ
\label{eqGU}
S_{ij}=<e_i^N,e_j^N>_N
\EE
Let $\sigma_{min}$ be the smallest eigenvalue of $G$ relative to $S$, i.e.
$$G\xi = \sigma_{min}S\xi$$
for some nonzero $\xi \in \Real^s$. Using Lagrange multipliers we can prove 
$$(\epsilon^N)^2=\sigma_{min}\rho$$
Therefore, (\ref{unobsgramian}) still holds true. 

For the heat equation approximated using (\ref{eqheatode}), the associated mappings can be defined by
$$\begin{array}{lll}
P^N(u)=\MT u^N_1,u^N_2,\cdots,u^N_N\EM^T, &u^N_k=\Fr{2}{L}\ds\int_0^{2\pi} u(x)\sin \left(\Fr{k\pi x}{L}\right)dx\\
\Phi^N(u^N)=\ds\sum_{k=1}^{N} u_k^N\sin \left(\Fr{k\pi x}{L}\right)
\end{array}$$
If we want to find the observability of the first $s$ modes, Definition \ref{def2} is equivalent to the analysis using the traditional observability gramian for $N=s$. In fact, for all $N\geq s$, $G$ is a constant matrix and 
$$G=\ds\int_0^Te^{(A^s)'t}(C^s)'C^se^{A^st}dt$$
Therefore, $\epsilon^N=\epsilon^s$ for all $N\geq s$ and $\epsilon^N$ is consistent. 

The idea of gramian matrix can be applied to nonlinear systems as a first order approximation of observability. Let $u^+_i(t)$ and $u^-_i(t)$ be solutions of the ODE in (\ref{hatuN}) with initial value
$$u^{\pm}_i(0)=u^N(0)\pm \rho e^N_i$$ 
If the system is linear, we have 
$$\begin{array}{lll}
\calH\circ\Phi^N(e^{-A^Nt}e^N_i)\\
=\Fr{1}{2\rho}\left(\calH\circ\Phi^N(u^+_i(t)) - \calH\circ\Phi^N(u^-_i(t))\right)
\end{array}$$
Therefore, if the ODE in (\ref{hatuN}) is nonlinear, we can define 
$$G_{ij}=\Fr{1}{4\rho^2}<\calH\circ\Phi^N(u^+_i(t)) - \calH\circ\Phi^N(u^-_i(t)), \calH\circ\Phi^N(u^+_j(t)) - \calH\circ\Phi^N(u^-_j(t))>_N$$
Let $\sigma_{min}$ be the smallest eigenvalue of $G$ relative to $S$, then (\ref{unobsgramian}) is the first order approximation of the unobservability index. In this case, $G$ is called an
{\it empirical gramian matrix}. This approach is inspired by the computational method in \cite{kreneride}.  

\section{Nonlinear systems}
\label{sec6}
\setcounter{equation}{0}

Consider a nonlinear initial value problem
\EQ
\label{eqpdemodelnon}
u_t=F(t,u,u_x,\cdots), &\mbox{in } \Omega \times (t_0, t_1]\\
\calB u=0 &\mbox{on } \partial \Omega_b \times (t_0, t_1]\\
u=u_0 & \mbox{in } \Omega \times \{ t=0\}\\
y(t)=\calH (u(\cdot,t))
\EE
where $\Omega$ is an open set in $\Real^n$, $F$ is a continuous function of $t$, $u$, and its derivatives with respect to $x$. The function $u(\cdot,t)$ belongs to $D$, which is a subspace of a Banach space $X$. The boundary condition is defined by a set of linear boundary operators on a part of $\partial \Omega$, denoted by $\partial \Omega_b$. To take into consideration of the boundary condition, we define
$$D_B=\{v\in D|\; \calB v =0 \mbox{ on } \partial \Omega_b\}$$
In the following, $u(t)$ represents $u(\cdot,t)$ in $D_B$. A solution $u(t)$ for this problem is a $\L2O$ valued function that is continuous in $[t_0,t_1]$, $du/dt$ exists and is continuous for all $t\in (t_0,t_1]$, satisfying $u(0)=u_0$, and $u(t)\in \DB$ for all $t\in (0,t_1]$. We assume that (\ref{eqpdemodel}) is  well posed. Proving the well-posedness of nonlinear PDEs is not easy. Nevertheless, for well posed problems the consistency of observability is guaranteed. In (\ref{eqpdemodel}), $y(t)=\calH (u(\cdot,t))$ represents the output of the system in which $\calH$ is a mapping, linear or nonlinear, from $\L2O$ to $\Real^p$. The output $y(t)\in C([t_0,t_1])$ has a norm denoted by $|| y||_Y$. Suppose 
$$W=\mbox{span}\{ e_1,e_2,\cdots,e_s\}$$
is the space for estimation, which is a finite dimensional subspace of $D_B$ generated by a basis $\{ e_1,e_2,\cdots,e_s\}$. For any solution, $u$, of (\ref{eqpdemodelnon}), similar to (\ref{decomposition}) we can define its best estimate $u_w\in W$ and its remainder $u_r$. Following Definition \ref{def1}, the ratio $\rho/\epsilon$ represents the unobservability index of $u_w$, where
\EQ
\label{eqcostnon}
\epsilon = \inf ||\calH (\hat u)-\calH (u_w)||_Y
\EE
and $\hat u$ satisfies
\EQ
\label{hatunon}
\hat u_t=F(t,u,u_x,\cdots)\\
\hat u (t)\in \DB\mbox{ for all } t\in (t_0,t_1]\\
\hat u(0)\in W\\
||\hat u(0) -u_w(0)||_X=\rho 
\EE

The approximation scheme consists of a sequence of ODEs
\EQ
\label{eqodemodelnon}
du^N/dt=F^N(t,u^N), & u^N\in \Real^N\\
u^N(0)=u^N_0
\EE
and two sequences of linear mappings
\EQ
\label{mappingnon}
P^N: D \rightarrow \Real^N\\
\Phi^N: \Real^N \rightarrow \L2O
\EE
The norm in $\Real^N$ is represented by $||\cdot ||_N$. Similar to the previous section, we assume that the approximation scheme is  well posed. In each $\Real^N$, the space for estimation is defined by 
$$\begin{array}{lll}
W^N=P^N(W)\\
=\mbox{span}\{ e_1^N, e_2^N,\cdots,e_s^N\}
\end{array}$$
For any trajectory $u^N(t)$ of (\ref{eqodemodelnon}) with $u^N(0)\in W^N$, its unobservability index is defined by $\rho/\epsilon^N$ in which 
$$\epsilon^N = \inf ||\calH \circ\Phi^N(\hat u^N)-\calH \circ\Phi^N(u^N)||_Y$$
subject to 
\EQ
\label{hatuNnon}
d\hat u^N/dt=F^N(t,u^N)\\
\hat u^N(0)\in W^N\\
||\hat u^N(0) -u^N(0)||_N=\rho 
\EE

For nonlinear systems, we need an additional assumption to guarantee the consistency. \\
\\
{\bf Convergence Assumption for $\Phi^N$}: Given any sequence 
$$\{ \Phi^N(v^N(t))\}_{N=N_0}^\infty$$ 
where $v^N(t)$ is a solution of (\ref{eqodemodelnon}) with $v^N(0)\in W^N$. If $\Phi^N(v^N(0))$ converges to $u(0)\in W_E$, where $u(t)$ is a solution of (\ref{eqpdemodelnon}), then $\Phi^N(v^N(t))$ converges to $u(t)$ uniformly for $t\in [t_0,t_1]$.\\

This assumption is stronger than the well-posedness assumption (\ref{mapping2}) because $v^N(0)$ is not required to equal $P^N( u(0))$. In many approximation schemes, the error satisfies 
$$||u(t)-\Phi^N(u^N(t))|| \leq \Fr{M}{N^\alpha}$$
for $u(0)$ from a bounded set in $W$, where $P^N(u(0))=u^N(0)$ and $\alpha >0$ is a rate of convergence. In this case, it can be proved that Convergence Assumption for $\Phi^N$ is satisfied.    

\begin{lemma}
\label{lemma2}
Suppose $\Phi^N$ satisfies Convergence Assumption. Given a sequence $\hat u^{N}(t)$, $N\geq N_0$, satisfying (\ref{hatuNnon}). There exists a subsequence, $\hat u^{N_k}(t)$, such that $\{ \Phi^{N_k}(\hat u^{N_k})\}_{k=1}^\infty$ converges uniformly to a solution of (\ref{hatunon}). 
\end{lemma}

{\it Proof }. For each $\hat u^N(t)$ in Lemma \ref{lemma2}, we have
$$\hat u^N(0)=\ds\sum_{i=1}^s a_i^N e_i^N$$
for some constants $a^N_1$, $a^N_2$, $\cdots$, $a^N_s$. 
From the condition on the initial value in (\ref{hatuNnon}), we know that the set 
$$\{||\ds\sum_{i=1}^s a_i^N e_i^N||_N, \;\; N\geq N_0\}$$
is bounded. Then it can be proved that the sequence $\{(a_1^N,a_2^N,\cdots,a_s^N)^T\}_{N=N_0}^\infty$ has a bounded subsequence that converges, 
$$\{(a_1^{N_k},a_2^{N_k},\cdots,a_s^{N_k})^T\}_{k=1}^\infty\rightarrow (a_1,a_2,\cdots,a_s)^T$$ 
Therefore
$$\begin{array}{rcl}
&&\ds\lim_{k\rightarrow \infty} \Phi^{N_k}(\hat u^{N_k}(0))\\
&=&\ds\lim_{k\rightarrow \infty} \ds\sum_{i=1}^s a_i^{N_k} \Phi^{N_k}(e_i^{N_k}))\\
&=&\ds\sum_{i=1}^s a_ie_i\\
&\triangleq& \hat u(0)
\end{array}$$
From Convergence Assumption for $\Phi^N$, $\{ \Phi^{N_k}(\hat u^{N_k})\}_{k=1}^\infty$ converges to $\hat u(t)$, a solution of the PDE in (\ref{hatunon}). 
The limit converges uniformly. Because 
$$\begin{array}{lll}
||\hat u^{N_k}(0)-u^{N_k}(0)||_{N_k}
=\rho
\end{array}$$ 
and because of the consistency of the norms, we have 
$$\begin{array}{lll}
|| \hat u(0)-u(0)||_X\\
=\ds\lim_{k\rightarrow \infty} || \ds\sum_{i=1}^s a_i e_i^{N_k} -u^{N_k}(0)||_{N_k}\\
=\ds\lim_{k\rightarrow \infty} || \ds\sum_{i=1}^s a_i^{N_k} e_i^{N_k} -u^{N_k}(0)||_{N_k}\\
=\rho
\end{array}$$ 
Therefore, $\{ \Phi^{N_k}(\hat u^{N_k})\}_{k=1}^\infty$ 
converges uniformly to a solution of (\ref{hatunon}).
$\square$\\

\begin{theorem} 
\label{theorem2}
Suppose the initial value problem (\ref{eqpdemodelnon}) and its approximation scheme (\ref{eqodemodelnon})-(\ref{mappingnon}) are well posed. Suppose $\Phi^N$ satisfies Convergence Assumption. Suppose $\calH$ satisfies Output Continuity Assumption. Then, 
\EQ
\label{consistnon}
\ds\lim_{N\rightarrow \infty} \epsilon^N =\epsilon
\EE
\end{theorem}

{\it Proof}. Most part of the proof is similar to that of Theorem \ref{theorem1} except for a few places where the linearity property is used. First, we prove
\EQ
\label{lowlimnon}
\liminf \epsilon^N \geq \epsilon
\EE
Suppose this is not true, then $\liminf \epsilon^N < \epsilon$. There exists $\alpha >0$ and a subsequence $N_k\rightarrow \infty$ so that
$$\epsilon^{N_k} < \epsilon -\alpha$$
for all $N_k$. From the definition of $\epsilon^N$, there exist $\hat u^{N_k}(t)$ satisfying (\ref{hatuNnon}) such that 
$$||\calH\circ\Phi^{N_k}(\hat u^{N_k})-\calH\circ\Phi^{N_k}(u^{N_k})||_Y < \epsilon - \alpha$$
From Lemma \ref{lemma2}, we can assume that $\Phi^{N_k}(u^{N_k}(t))$ converges to $\hat u(t)$ uniformly, where $\hat u(t)$ satisfies (\ref{hatunon}). From Output Continuity Assumption, 
$$\lim_{k\rightarrow \infty} ||\calH(\Phi^{N_k}(\hat u^{N_k})-\calH(\Phi^{N_k}(u^{N_k})||_Y =||\calH(\hat u)-\calH(u)||_Y \leq \epsilon - \alpha$$
However, from the definition of $\epsilon$, we know 
$$ \epsilon \leq ||\calH(\hat u)-\calH(u)||_Y$$
A contradiction is found. Therefore, (\ref{lowlimnon}) must hold. 

In the next step, we prove
\EQ
\label{uplimnon}
\limsup \epsilon^N \leq \epsilon
\EE
It is adequate to prove the following statement: for any $\alpha>0$, there exists $N_1>0$ so that 
\EQ
\label{ineqnon}
\epsilon^N < \epsilon +\alpha
\EE
for all $N\geq N_1$. From the definition of $\epsilon$, there exists $\hat u$ satisfying (\ref{hatunon}) so that 
\EQ
\label{eqad1non}
||\calH(\hat u)-\calH(u)||_Y < \epsilon +\alpha
\EE
Let $\hat u^N$ be a solution of the ODE
$$du^N/dt=F^N(t,u^N)$$
with an initial value
$$\hat u^N(0)=P^N(\hat u(0))$$ 
Then the following limit converges uniformly
\EQ
\label{eqlimnon}
\ds\lim_{N\rightarrow \infty} ||\Phi^N(\hat u^N(t))-\hat u(t)||_X = 0
\EE
A problem with $\hat u^N(0)$ is that its distance to $u^N(0)$ may not be $\rho$, which is required by (\ref{hatuNnon}). Let $\bar u(t)$ be a solution of (\ref{eqodemodelnon}) with an initial value 
$$\begin{array}{lll}
\bar u^N(0)=\gamma_N(\hat u^N(0)-u^N(0))+u^N(0)\\
\gamma_N=\Fr{\rho}{||\hat u^N(0)-u^N(0)||_N}
\end{array}$$
Obviously 
$$||\bar u^N(0)-u^N(0)||_N=\rho$$
and $\bar u(t)$ satisfies (\ref{hatuNnon}). 
Due to the consistency of the norms and the fact $||\hat u(0)-u(0)||_X=\rho$, we know
\EQ
\label{limgamma}
\lim_{N\rightarrow \infty} \gamma_N=1
\EE
From (\ref{eqlimnon}) and (\ref{limgamma}), we have
$$\begin{array}{lll}
\ds\lim_{N\rightarrow \infty} \Phi^N(\bar u^N(0))\\
=\ds\lim_{N\rightarrow \infty}\left(\gamma_N(\Phi^N(\hat u^N(0))-\Phi^N(u^N(0)))+\Phi^N(u^N(0))\right)\\
=\hat u(0)
\end{array}$$
By Convergence Assumption for $\Phi^N$, the limit 
$$\ds\lim_{N\rightarrow \infty} \Phi^N(\bar u^N(t))=\hat u(t)$$
converges uniformly on the interval $t \in [t_0,t_1]$. 
Then Output Continuity Assumption implies
$$\begin{array}{lll}
\ds\lim_{N\rightarrow \infty} ||\calH(\Phi^N(\bar u^N))-\calH(\Phi^N(u^N))||_Y\\
=||\calH(\hat u)-\calH (u)||_Y\\
<\epsilon + \alpha
\end{array}$$
There exits $N_1>0$ so that
$$||\calH(\Phi^N(\bar u^N))-\calH(\Phi^N(u^N))||_Y < \epsilon+\alpha$$
for all $N\geq N_1$. From the definition of $\epsilon^N$, we know
$$\epsilon^N \leq ||\calH(\Phi^N(\bar u^N))-\calH(\Phi^N(u^N))||_Y< \epsilon + \alpha$$
for all $N>N_1$. Therefore, (\ref{uplimnon}) holds. 

To summarize, the inequalities (\ref{lowlimnon}) and (\ref{uplimnon}) imply 
$$\lim_{N\rightarrow \infty} \epsilon^N =\epsilon$$
$\square$

\section{Example}
\label{sec7}
\begin{figure}[!ht]
	\begin{center}
		\includegraphics[width=3.0in,height=3.0in]{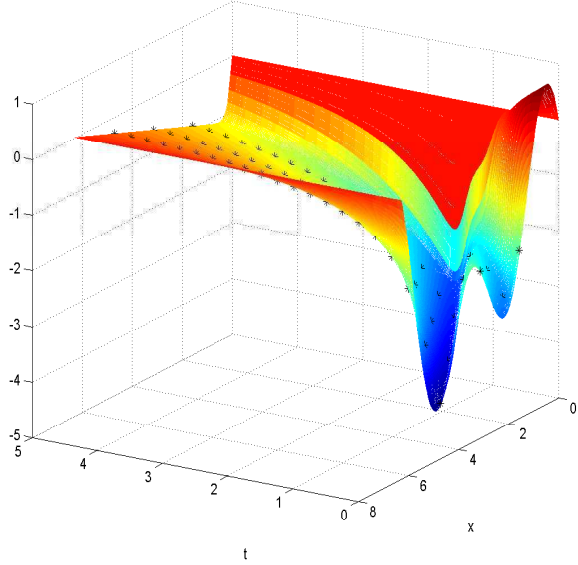} 
		\caption{A solution of Burgers' equation with sensor measurements}
		\label{Figburgers}
		\end{center}
\end{figure}

\setcounter{equation}{0}
Consider the following Burgers' equation
\EQ
\label{pde}
\Fr{\partial u(x,t)}{\partial t}+u(x,t)\Fr{\partial u(x,t)}{\partial x}=\kappa \Fr{\partial^2 u(x,t)}{\partial x^2} \\
u(x,0)=u_0(x), & x\in [0, L]\\
u(0,t)=u(L,t)=0, & t\in [0,T]
\EE
where $L=2\pi$, $T=5$, and $\kappa=0.14$. The following output represents three sensors that measure the value of $u(x,t)$ at fixed locations 
\EQ
\label{pdeoutput}
\MT y_1(t_k) \\ y_2(t_k)\\ y_3(t_k)\EM=\MT u(\frac{L}{4},t_k)\\ u(\frac{2L}{4},t_k)\\ u(\frac{3L}{4},t_k)\EM, & t_k=k\Delta t, \;\; k=0,1,2,\cdots,N_t
\EE
where $\Delta t=T/N_t$, $N_t=20$. Figure \ref{Figburgers} shows a solution with discrete sensor measurements marked by the stars. In the output space
$$||y||_Y=\left(\dsum_{k=0}^{N_t}(y_1^2(k)+y_2^2(k)+y_3^2(k))\right)^{1/2}$$

The approximation scheme is based on equally spaced grid-points
$$x_0=0 < x_1 <\cdots < x_N=L,$$ 
where 
$$\Delta x=x_{i+1}-x_i=L/N.$$  
System (\ref{pde}) is discretized using a central difference method
\EQ
\label{odemodel}
\dot u_1^N=-u_1^N\Fr{u_2^N-u_0^N}{2\Delta x}+\kappa \Fr{u_{2}^N-2u_1^N}{\Delta x^2}\\
\dot u_2^N=-u_2^N\Fr{u_3^N-u_1^N}{2\Delta x}+\kappa \Fr{u_{3}^N+u_{1}^N-2u_2^N}{\Delta x^2}\\
\;\;\;\;\;\vdots\\
\dot u_{N-1}^N=-u_{N-1}^N\Fr{u^N_{N}-u_{N-2}^N}{2\Delta x}+\kappa \Fr{u_N^N+u_{N-2}^N-2u_{N-1}^N}{\Delta x^2}
\EE
where $u^N_0=u^N_{N}=0$. For any $v(x)\in C^1([0,L])$, we define
$$P^N(v)=\MT v(x_1) & v(x_2) & \cdots v(x_{N-1})\EM \in \Real^{N-1}$$
For any $v^N\in \Real^{N-1}$, define
$$\Phi^N(v^N)=\frac{a_0}{2}+\dsum_{k=1}^{N/2-1} \left(a_k\cos\left(\frac{2\pi k}{L}x\right)+b_k\sin\left(\frac{2\pi k}{L}x\right)\right)+\frac{a_{N/2}}{2}\cos\left(\frac{\pi N}{L}x\right)$$
where
\EQ
\label{eqfouriercoef}
a_k=\Fr{2}{N}\dsum_{j=1}^{N-1} v^N_j\cos\left(\frac{2\pi k}{L}x_j\right), & k=0,1,2,\cdots, N/2\\
b_k=\Fr{2}{N}\dsum_{j=1}^{N-1} v^N_j\sin\left(\frac{2\pi k}{L}x_j\right), & k=1,2,\cdots, N/2-1
\EE
We adopt $L^2$-norm in $C^1[0, L]$. For any vector $v^N \in \Real^{N-1}$, we defined its norm using the Fourier coefficients (\ref{eqfouriercoef}),
$$||v^N||_N=(a_0^2/2+\dsum_{k=1}^{N/2-1} (a_k^2+b_k^2)+a_{N/2}^2/4)\pi$$
From the orthogonality of trigonometic functions,
$$||\Phi^N(v^N)||_{L^2} =  ||v^N||_N$$

The space for estimation is defined to be 
$$W=\left\{\alpha_0/2+ \dsum_{k=1}^{K_F} \left(\alpha_k\cos(\Fr{2k\pi}{L}x)+\beta_k \sin(\Fr{2k\pi}{L}x)\right)\left |\begin{array}{lll}\alpha_k, \beta_k \in \Real\\ \alpha_0/2+ \dsum_{k=1}^{K_F} \alpha_k=0\end{array}\right. \right\} $$
In this section, $K_F=2$.  This means that we want to find the observability for the first five modes in the Fourier expansion of $u(0)$. Or equivalently, we would like to find the observability of 
$$\MT\alpha_0 & \alpha_1 &\beta_1 &\alpha_2 &\beta_2\EM$$
Define 
$$X^N=\MT x_1 & x_2&\cdots &x_{N-1}\EM^T$$
then
$$W^N=\left\{\alpha_0/2+ \dsum_{k=1}^{K_F} \left(\alpha_k\cos(\Fr{2k\pi}{L}X^N)+\beta_k \sin(\Fr{2k\pi}{L}X^N)\right)\left |\;  \begin{array}{lll}\alpha_k, \beta_k \in \Real\\ \alpha_0/2+ \dsum_{k=1}^{K_F} \alpha_k=0\end{array}\right. \right\} $$

\begin{figure}[!ht]
	\begin{center}
		\includegraphics[width=3.0in,height=2.0in]{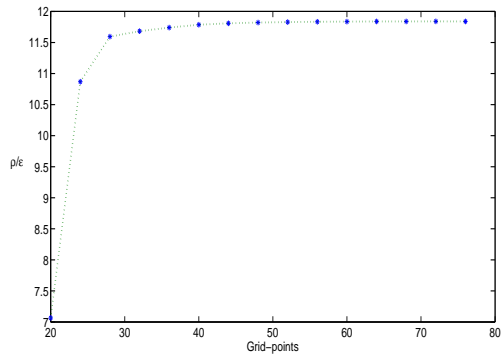} 
		\caption{Consistency for Burgers' equation}
		\label{figburger}
		\end{center}
\end{figure}

For this simulation, the nominal trajectory has the following initial value 
$$u_0(x)=-2+\cos(x)+\sin(x)+\cos(2x)+\sin(2x)$$
Its solution is shown in Figure \ref{Figburgers}. To approximate its observability, we apply the empirical gramian method to (\ref{odemodel}) in the space $W^N$. The consistency of observability is verified by the results. The ratio $\rho/\epsilon^N$ is approximated for $N=4k$, $5\leq k \leq 19$. The value of unobservability index approaches (Figure \ref{figburger}) $$\rho/\epsilon=11.83$$  

\section{Conclusions}
A definition of observability using dynamic optimization is introduced for PDEs. The advantage of this definition is to resolve several issues and concerns about observability in a unified framework. More specifically, using the concept one can achieve a quantitative measure of partial observability for PDEs. Furthermore, the observability can be numerically approximated. A practical feature of this definition for infinite dimensional systems is that the observability can be numerically computed using well posed approximation schemes. It is mathematically proved that the approximated observability is consistent with the observability of the original PDE. A first order approximation is derived using empirical gramian matrices. The consistency is verified using an example of a Burgers' equation.  \\
\\
{\bf Acknowledgement}:  The author would like to express his gratitude to Professor Arthur J. Krener for his comments and discussions on the observability gramian and unobservability index. The author would like to thank Dr. Liang Xu for his suggestions and comments on potential applications in numerical weather prediction.

\end{document}